# Markov Boundary Discovery with Ridge Regularized Linear Models

**Eric V. Strobl**  EVS17@PITT.EDU
**Shyam Visweswaran**  SHV3@PITT.EDU
*Center for Causal Discovery*
*Department of Biomedical Informatics*
*University of Pittsburgh School of Medicine*
*5607 Baum Boulevard*
*Pittsburgh, PA 15206, USA*

**Abstract**

Ridge regularized linear models (RRLMs), such as ridge regression and the SVM, are a popular group of methods that are used in conjunction with coefficient hypothesis testing to discover explanatory variables with a significant multivariate association to a response. However, many investigators are reluctant to draw causal interpretations of the selected variables due to the incomplete knowledge of the capabilities of RRLMs in causal inference. Under reasonable assumptions, we show that a modified form of RRLMs can get "very close" to identifying a subset of the Markov boundary by providing a worst-case bound on the space of possible solutions. The results hold for any convex loss, even when the underlying functional relationship is nonlinear, and the solution is not unique. Our approach combines ideas in Markov boundary and sufficient dimension reduction theory. Experimental results show that the modified RRLMs are competitive against state-of-the-art algorithms in discovering part of the Markov boundary from gene expression data.

**Keywords:** Markov boundary, ridge regularization, linear models

## 1. Introduction

Inferring causal relationships between random variables from observational data is a difficult task. However, algorithms that can accurately identify causal relations may significantly speed up scientific investigations. In many applications, scientists are primarily interested in recovering a restricted group of variables which share causal relationships with a given response variable. For example, biologists may be interested in using gene expression levels to identify new targets which effect cancer survivability. Algorithms which can discover variables that share causal relationships with the response can thus be very useful in practice.

Under mild assumptions (details in Section 2.3), the Markov boundary is a unique group of random variables located within the close vicinity of the response variable in a casual directed graph, where causal relations are modeled by directed edges between variables. In a causal directed acyclic graph, the Markov boundary consists of the response's direct causes, direct effects, and direct causes of the direct effects (Neapolitan 2004, Tsamardinos and Aliferis 2003). Otherwise, the Markov boundary can be a subset of those variables (Neapolitan 2004, Pearl 1988). Note that certain manipulations of the direct causes may modify the distribution of the response, and certain manipulations of the direct effects or direct causes of the direct effects may modify the distributions of the response's consequences. The Markov boundary can thus be a small set of variables that share important causal relationships with a variable of interest.

Variables in the Markov boundary exhibit the following conditional independence relation with other variables in the dataset:

$$Y \perp\!\!\!\perp (X \setminus M) | M, \qquad (E1)$$

where $Y$ denotes the response variable, $X$ all variables other than $Y$, $M$ the Markov boundary, and $\perp\!\!\!\perp$ statistical independence. Unlike the Markov blanket, the Markov boundary is defined so that no proper subset of $M$ also satisfies (E1) (Statnikov et al. 2013). This fact seems to suggest that Markov boundary discovery is closely related to the variable selection problem. This intuition is indeed true (details in Sections 2.4 and 3.2), but most variable selection methods unfortunately do not attempt to recover the Markov boundary. Their primary goal is to identify a set of variables that maximizes the predictive accuracy of a model, so the causal interpretation of the selected variables has often been heuristic in nature. This is unfortunate because there is often a strong desire to interpret the selected variables as having some sort of causal relationship with the response. Whole research programs have in fact been built on experiments motivated by causal hypotheses that were generated by non-causal variable selection algorithms (Zhou, Kao, and Wong 2002, Li et al. 2001, Holmes, Durbin, and Winston 2000, Eisen et al. 1998). It is thus important to develop a theoretical understanding of both the advantages and disadvantages of using particular selection strategies for causal inference.

Ridge regularized linear models (RRLMs), such as ridge regression and the support vector machine (SVM)[1], are popular methods that are used to identify variables with a multivariate association with the response. To understand ridge or $L^2$-regularization intuitively, recall that the mean squared error (MSE) of the estimator $\hat{\beta}$ with respect to the variable coefficients $\beta \in \mathbb{R}^p$ can be decomposed into a bias term and a variance term. Ridge regularization penalizes the size of the variable coefficients by shrinking the values of the estimate $\hat{\beta}$ towards zero; this introduces bias but reduces the variance of the estimate. The method is thus particularly useful in the high dimensional low sample size (HDLSS) context, where the variance of $\hat{\beta}$ can be large. In fact, in ridge regression, there always exists a regularization constant $\lambda \in \mathbb{R}^+$ such that the MSE of $\hat{\beta}$ estimated via ridge regression is smaller than the MSE of $\hat{\beta}$ estimated via ordinary least squares (Hoerl and Kennard 1970). Moreover, in the context of linear SVMs, recall that the distance between the support vectors is $2/\|\beta\|_2$ in the separable binary case. As a result, minimizing $\|\beta\|_2$ is equivalent to maximizing the margin distance, an important quantity which can be used to bound the SVM's generalization error (Vapnik and Chapelle 2000).

RRLMs are also popular for several other reasons besides variance reduction and margin maximization. First, most RRLMs can be optimized efficiently (e.g., Glmnet (Qian et al. 2013) and LIBSVM (Chang and Lin 2011)). As a result, permutation analyses can be employed to recover p-values for each variable; these p-values are often the primary quantity of interest in both variable selection and ranking, since they *may* help quantify the probability that subsequently conducted experiments will fail to identify an effect. RRLMs are also used in the HDLSS context, where detection of nonlinearities can be difficult, and linear models are known to often predict better than or comparably with nonlinear methods. Finally, variable importance can be easily assessed in linear models using the absolute values of the coefficients, whereas complex strategies may be needed to recover variable importance values from nonlinear methods. This implies that backward elimination can also be carried out efficiently with RRLMs by iteratively eliminating the variable associated with the coefficient with the smallest absolute value. Together, these strengths highlight the importance of developing an understanding of the causal inference capabilities of RRLMs.

---

[1] The linear SVM can be implemented with $L^1$ and/or $L^2$-regularization. In this paper, we refer to the original and most popular version of the SVM with $L^2$-regularization and the hinge loss (Wu and Liu 2007, Chapelle , Cortes and Vapnik 1995).

In this paper, we characterize the ability of RRLMs in identifying the Markov boundary. Our results show that a modified form of RRLMs can get "very close" to identifying a subset of the Markov boundary by providing a worst-case bound on the space of possible solutions. This result holds even if the underlying functional relationship between the explanatory variables and response variable is nonlinear. Furthermore, the solution to each RRLM's optimization problem does not need to be unique. The theory is based on the fact that the expectation of the explanatory variables can almost always be written as a linear function of their lower dimensional projections $\eta^T X$ in the high dimensional context (Hall and Li 1993). We also present experimental results which demonstrate that the modified RRLMs are competitive when compared with more sophisticated methods for discovering a subset of the Markov boundary from HDLSS gene expression data.

## 2 Background

### 2.1 Notation

Unless specified otherwise, upper case letters in italics denote random variables (e.g., $A$, $B$, $C$) and upper case bold letters in italics denote random variable sets or random column vectors ($\boldsymbol{A}$, $\boldsymbol{B}$, $\boldsymbol{C}$). We reserve $Y$ for the response variable, and $\boldsymbol{X}$ for the set or column vector of all $p$ explanatory variables. We also use the notation $\boldsymbol{A} \subseteq \boldsymbol{B}$ to mean $\boldsymbol{A}$ is a subset of $\boldsymbol{B}$, and $\boldsymbol{A} \subset \boldsymbol{B}$ to mean $\boldsymbol{A}$ is a proper subset of $\boldsymbol{B}$.

We will be arranging both vectors and matrices in this paper, so we use the following notation to distinguish between operations involving rows versus columns. If $\boldsymbol{A}$ is a $a \times 1$ vector and $\boldsymbol{B}$ is a $b \times 1$ vector, then we use the notation $(\boldsymbol{A}; \boldsymbol{B})$ to denote that $\boldsymbol{B}$ is appended to the rows of $\boldsymbol{A}$ to create a $(a + b) \times 1$ column vector. On the other hand, we use the notation $(\boldsymbol{A}, \boldsymbol{B})$ to denote that $\boldsymbol{B}$ is appended to the columns of $\boldsymbol{A}$ to create a $a \times 2$ matrix when $a = b$.

A graph consists of nodes and edges, where nodes are random variables and directed edges are causal relationships between two variables. If a graph contains a directed edge between two variables $A$ and $B$ such that $A \to B$, then $A$ is a direct cause or parent of $B$, and $B$ is a direct effect or child of $A$. A node $A$ is a direct cause of the direct effect or a spouse of $B$, if both $A$ and $B$ share a common child node.

We write $\boldsymbol{A} \perp\!\!\!\perp \boldsymbol{B} | \boldsymbol{C}$ when two sets of variables $\boldsymbol{A}$ and $\boldsymbol{B}$ are conditionally independent given a third set of variables $\boldsymbol{C}$. We also write $\boldsymbol{A} \perp\!\!\!\perp \boldsymbol{B}$ when the conditioning set is empty to denote that $\boldsymbol{A}$ and $\boldsymbol{B}$ are independent.

### 2.2 Distribution Theory

The following is an important property of a probability distribution, which we will often refer to in this paper:

**Definition 2.2.1. (Intersection Property)** Let $\boldsymbol{A}$, $\boldsymbol{B}$, $\boldsymbol{C}$ and $\boldsymbol{D}$ be any four subsets of potentially overlapping variables from $(Y; \boldsymbol{X})$ with joint probability distribution $\mathbb{P}$. The distribution $\mathbb{P}$ is said to satisfy the intersection property if $\boldsymbol{A} \perp\!\!\!\perp \boldsymbol{B}|(\boldsymbol{C} \cup \boldsymbol{D})$ and $\boldsymbol{A} \perp\!\!\!\perp \boldsymbol{D}|(\boldsymbol{C} \cup \boldsymbol{B}) \implies \boldsymbol{A} \perp\!\!\!\perp (\boldsymbol{B} \cup \boldsymbol{D})|\boldsymbol{C}$. ∎

There are two sufficient conditions for the intersection property:

**Proposition 2.2.2.** If $\mathbb{P}$ has a strictly positive density or has a continuous density with path-connected support, then $\mathbb{P}$ satisfies the intersection property. ∎

By continuous density with path-connected support, we mean that an absolutely continuous distribution has a continuous density whose support can be connected by axis-parallel lines; this is in fact both a necessary and sufficient condition. Note that we will always assume that a distribution has a density in this paper, so we will use the words "distribution" and "density" interchangeably. For the proof of Proposition 2.2.2, the part regarding strict positivity is provided in (Pearl 1988, Dawid 1979), and the part regarding path-connected support is provided in (Peters 2014).

We believe that these sufficient conditions are mild and are satisfied in most distributions encountered in practice. Strict positivity is considered reasonable whenever there is uncertainty about the data (Pearl 1988). Moreover, it is reasonable to assume that almost all distributions with continuous densities have path-connected support in biology. We nonetheless must acknowledge that several investigators have challenged these assumptions (e.g., (Statnikov et al. 2013, Lemeire and Dominik 2010, Lemeire, Meganck, and Cartella 2010)). We do not claim that strict positivity or path-connected support holds in all biological datasets; however, we believe that the empirical results used to refute these assumptions are still inconclusive, since some of the same investigators studying HDLSS molecular data have warned that their experimental results can also hold when sample sizes are insufficient (Statnikov et al. 2013, Dougherty and Brun 2006, Ein-Dor, Zuk, and Domany 2006, Statnikov and Aliferis 2010).

## 2.3 Markov Boundary Theory

We now define a special group of variables which helps characterize a conditional independence relation concerning the response variable $Y$:

**Definition 2.3.1: Markov blanket.** A Markov blanket $M$ of a response variable $Y$ in the joint probability distribution $\mathbb{P}$ over variables $X$ is a set of variables conditioned on which all other variables are independent of $Y$; that is, for every $F \subseteq (X \backslash M), Y \perp\!\!\!\perp F | M$. ∎

The Markov blanket of $Y$ is shown to exhibit the following relation using (E1):

**Proposition 2.3.2.** We have $Y \perp\!\!\!\perp (X \backslash M) | M \Leftrightarrow Y \perp\!\!\!\perp X | M$. ∎

The proof follows directly from Proposition 4.4 of (Cook 1998). Note that a trivial Markov blanket of $Y$ is $X$. As a result, we more specifically define a minimal Markov blanket:

**Definition 2.3.3: Markov boundary.** If no proper subset of $M$ satisfies the definition of a Markov blanket of $Y$, then $M$ is the Markov boundary of $Y$. ∎

We now have a sufficient condition for the uniqueness of the Markov boundary of $Y$:

**Theorem 2.3.4.** (Pearl 1988) If a joint probability distribution $\mathbb{P}$ over variables $(Y; X)$ satisfies the intersection property, then for each $V \in (Y; X)$, there exists a unique Markov boundary of $V$. ∎

Note that some investigators have studied violations of the intersection property, such as in deterministic relations, where uniqueness of the Markov boundary may not hold (Statnikov et al. 2013, Lemeire and Dominik 2010). However, we assume that the intersection property holds in this paper.

We define an important condition which we will use in the next theorem:

**Definition 2.3.5: Global Markov condition.** (Richardson and Spirtes 1999, 2002). The joint probability distribution $\mathbb{P}$ over variables $H$ satisfies the global Markov condition for a directed graph $\mathbb{G} = <H, \mathbb{E}>$ if and only if for any three disjoint subsets of variables $A$, $B$ and $C$ from $H$, if $A$ is d-separated from $B$ given $C$ in $\mathbb{G}$ then $A$ is conditionally independent of $B$ given $C$ in $\mathbb{P}$. ∎

The following theorem guarantees that the Markov boundary is a subset of the children, parents and spouses of the response variable, if the global Markov condition is satisfied:

**Theorem 2.3.6.** (Neapolitan 2004, Pearl 1988) If a joint probability distribution $\mathbb{P}$ satisfies the global Markov condition for directed graph $\mathbb{G}$, then the set of children, parents and spouses of $Y$ is a Markov blanket of $Y$. ∎

## 2.4   The Markov Blanket and Variable Selection

This section connects the ideas of a Markov blanket to generic variable selection. We first consider an optimal set of variables for predicting $Y$.

**Definition 2.4.1: Optimal predictor.** (Statnikov et al. 2013) Given a dataset $\mathbb{D}$ (a sample from distribution $\mathbb{P}$) for variables $(Y; \boldsymbol{X})$, a learning algorithm $\mathbb{L}$, and a performance metric $\mathbb{M}$ to assess the learner's models, a variable set $\boldsymbol{V} \subseteq \boldsymbol{X}$ is an optimal predictor of $Y$ if $\boldsymbol{V}$ maximizes $\mathbb{M}$ for predicting $Y$ using learner $\mathbb{L}$ in the dataset $\mathbb{D}$. ∎

Examples of performance metrics include the log-likelihood, negative mean squared error, and negative hinge loss. We can equivalently consider minimizing a loss functional instead of maximizing a performance metric. The following theorem provides necessary and sufficient conditions for the Markov blanket to be an optimal predictor.

**Theorem 2.4.2.** (Statnikov et al. 2013) If $\mathbb{M}$ is a performance metric that is maximized only when the conditional probability distribution $\mathbb{P}(Y|\boldsymbol{X})$ is estimated accurately, and $\mathbb{L}$ is a learning algorithm that can approximate any conditional probability distribution, then $\boldsymbol{M}$ is a Markov blanket of $Y$ if and only if it is an optimal predictor of $Y$. ∎

## 2.5   Sufficient Dimension Reduction Theory

Sufficient dimension reduction (SDR) can be viewed as a generalization of Markov boundary discovery. SDR attempts to find a matrix $\eta \in \mathbb{R}^{p \times d}$, where $d \leq p$ such that:

$$Y \perp\!\!\!\perp \boldsymbol{X} | \eta^T \boldsymbol{X}. \tag{E2}$$

Notice that $\eta^T \boldsymbol{X}$ is a lower dimensional linear combination of $\boldsymbol{X}$. In contrast, Markov boundary discovery algorithms first consider $\eta = I_p$, and then use combinatorial optimization methods (e.g., constraint-based search) to set some of the diagonal elements to zero. We will see that the added flexibility of considering any lower dimensional linear combination of $\boldsymbol{X}$ which satisfies (E2) will allow us to convert the traditional combinatorial optimization approach of Markov boundary discovery to an easier, continuous optimization problem without sacrificing fidelity.

We provide the following two definitions from (Cook 1998). Note that we do not restrict ourselves to a continuous response or regression. We will use the notation $Y|\boldsymbol{X}$ to denote the response variable $Y$ given the explanatory variables $\boldsymbol{X}$; one can then speak of, for example, the probability distribution or expectation of $Y$ given $\boldsymbol{X}$, which we write as $\mathbb{P}(Y|\boldsymbol{X})$ and $E(Y|\boldsymbol{X})$, respectively. For instance, we may have:

$$Y|\boldsymbol{X} = X_1 + 2X_2 + \varepsilon,$$

where $\varepsilon$ is Gaussian noise and $\varepsilon \perp\!\!\!\perp \boldsymbol{X}$.

**Definition 2.5.1: Dimension reduction subspace.** If $Y \perp\!\!\!\perp X | \eta^T X$, then the column space of $\eta$ is called a dimension reduction subspace (DRS) for $Y|X$. ∎

The column space of $\eta$ refers to the space spanned by the columns of $\eta$ such that the columns form a basis. We will denote the column space of $\eta$ as $S(\eta)$. Also:

**Definition 2.5.2: Central dimension reduction subspace.** A subspace $S_{Y|X}$ is a central DRS for $Y|X$ if $S_{Y|X}$ is a DRS and $S_{Y|X} \subseteq S_{DRS}$ for all DRSs $S_{DRS}$. ∎

The following is an important sufficient condition for the existence of a central DRS:

**Theorem 2.5.3.** (Cook 1998) Let $S(\alpha)$ and $S(\phi)$ be DRSs for $Y|X$. If $X$ has a density $f(a) > 0$ for $a \in \Omega_X \subseteq \mathbb{R}^p$ and $f(a) = 0$ otherwise, and if $\Omega_X$ is a convex set, then $S(\alpha) \cap S(\phi)$ is a DRS. ∎

The central DRS for $Y|X$ thus exists if the density of $X$ has convex support because the intersection of all DRSs is the central DRS, $\cap S_{DRS} = S_{Y|X}$. A similar result also exists for discrete distributions with connected support (see Problem 6.5 in (Cook 1998)).

## 3 Theoretical Results

### 3.1 Summary

We provide a summary of the results of this section. Consider the following assumptions:

(1) The global Markov condition holds.
(2) The joint probability distribution of $(Y; X)$ satisfies the linear intersection property (defined in Section 3.3).
(3) $\Sigma_X$ is positive definite.
(4) $E(X|\eta^T X)$ is a linear function of $\eta^T X$ when $Y \perp\!\!\!\perp X | \eta^T X$.
(5) The matrix $\beta^*$ is a non-zero matrix, and the solution to the following optimization problem:

$$\underset{\alpha,\beta}{argmin}\ E\{u(\alpha + \beta^T X, Y)\} + \lambda tr(\beta^T \Sigma_X \beta).$$

where $\alpha \in \mathbb{R}^K$, $\beta \in \mathbb{R}^{p \times K}$, $u(\cdot,\cdot)$ is any convex functional, $tr(\cdot)$ denotes the trace, $\lambda \in \mathbb{R}^+$, and $K \in \mathbb{Z}^+$. Notice that a covariance matrix $\Sigma_X$ has been added into the ridge regularization, and the regularization is more specifically a squared Frobenius norm.

This section makes the following conclusions based on the above five assumptions. The linear intersection property (Assumption 2) ensures that a unique Markov boundary and a central DRS exist by Theorems 2.3.4 and 3.3.2, respectively. Theorem 3.5.2 then states that we can solve the optimization problem in (5) and find a $\beta^*$ such that $S(\beta^*) \subseteq S_{Y|X}$, since $\Sigma_X$ is positive definite (3), and $E(X|\eta^T X)$ is a linear function of $\eta^T X$ when $Y \perp\!\!\!\perp X | \eta^T X$ (4). Next, we can recover a nonempty subset of the Markov boundary from the non-zero coefficients of $\beta^*$ by Theorem 3.4.2, if we also assume that $\beta^*$ is a non-zero matrix (5). The global Markov condition (1) finally guarantees that the set of parents, children and spouses of $Y$ includes the Markov boundary of $Y$ by Theorem 2.3.6. As a result, if the above five assumptions hold, then $\beta^* \in S_{Y|X}$, and a nonempty subset of the parents, children and spouses of $Y$ is identifiable from $\beta^*$. Now, if Assumption 4 does not hold for all $\eta$ such that $Y \perp\!\!\!\perp X | \eta^T X$, then $S(\beta^*) \subseteq S(\eta)$ for those $\eta$ which do also by Theorem 3.5.2. Notice that the above argument does not require causal faithfulness but instead requires a variant of the intersection property. Moreover, the generality of the optimization problem in (5) implies that the result applies to RRLM-based regression, binary classification, and multi-class classification among possibly others.

### 3.2 The Markov Boundary and Variable Selection

We start our investigation by extending the theoretical connection between the Markov blanket and variable selection by more precisely considering the Markov boundary. Note that, like the Markov blanket of $Y$, a trivial optimal predictor is $X$. As a result, we can define a minimal and optimal predictor:

**Definition 3.2.1: Minimal and optimal predictor.** Let $V$ be an optimal predictor of $Y$. If no proper subset of $V$ satisfies the definition of an optimal predictor of $Y$, then $V$ is a minimal and optimal predictor of $Y$. ∎

Similar to Theorem 2.4.2, the following theorem provides necessary and sufficient conditions for a Markov boundary to be a minimal and optimal predictor.

**Theorem 3.2.2.** If $\mathbb{M}$ is a performance metric that is maximized only when the conditional probability distribution $\mathbb{P}(Y|X)$ is estimated accurately, and $\mathbb{L}$ is a learning algorithm that can approximate any conditional probability distribution, then $M$ is a Markov boundary of $Y$ if and only if it is a minimal and optimal predictor of $Y$.

**Proof.** First assume $M$ is a Markov boundary of $Y$. Then, $M$ is also a Markov blanket of $Y$, so $M$ is an optimal predictor of $Y$ by Theorem 2.4.2. Moreover, by the definition of a Markov boundary, no proper subset of $M$ is a Markov blanket of $Y$. As a result, no proper subset of $M$ also satisfies the definition of an optimal predictor of $Y$. We thus conclude that $M$ is a minimal and optimal predictor of $Y$.

The other direction follows similarly. Assume $G$ is a minimal and optimal predictor of $Y$. Then, $G$ is also an optimal predictor of $Y$, so $G$ is a Markov blanket of $Y$ by Theorem 2.4.2. Next, by the definition of a minimal and optimal predictor, no proper subset of $G$ is also an optimal predictor of $Y$. Thus, no proper subset of $G$ is also a Markov blanket of $Y$. We conclude that $G$ is a Markov boundary of $Y$. ∎

The above theorem implies that the Markov boundary of $Y$ is a solution to the variable selection problem when $\mathbb{P}(Y|X)$ needs to be accurately estimated. Unfortunately, this is not the case in many machine learning tasks. For example, in regression, the negative mean squared error is maximized when $E(Y|X)$ is accurately approximated – not the entire conditional distribution. Moreover, estimating the conditional expectation may require a subset of the variables in the Markov boundary of $Y$. As an example, let $X = (X_1, X_2)$ and $Y = f(X_1) + N(0, \rho(X_2)^2)$, where $f(\cdot)$ and $\rho(\cdot)$ are some fixed functions. Then, $\mathbb{P}(Y|X)$ depends on $X_1$ and $X_2$ but $E(Y|X)$ only depends on $X_1$. This example and Theorem 3.2.2 suggest that not all RRLM methods can identify the entire Markov boundary. We will see that this intuition is indeed true in Section 3.5.

We now restrict ourselves to tasks where $\mathbb{P}(Y|X)$ is estimated and define the following:

**Definition 3.2.3: Best predictor.** A best predictor is a minimal and optimal predictor with the smallest cardinality. ∎

Necessary and sufficient conditions regarding the Markov boundary and the best predictor can be similarly established by following the logic of the proof of Theorem 3.2.2:

**Theorem 3.2.4.** If $\mathbb{M}$ is a performance metric that is maximized only when $\mathbb{P}(Y|X)$ is estimated accurately, and $\mathbb{L}$ is a learning algorithm that can approximate any conditional probability distribution, then $M$ is a Markov boundary of $Y$ with smallest cardinality if and only if it is a best predictor of $Y$. ∎

If only one Markov boundary exists, such as when the intersection property holds, then the Markov boundary is the best predictor and vice versa. If multiple Markov boundaries exist, then the Markov boundary or boundaries with the smallest cardinality is/are the best predictor(s) and vice versa. Thus, the best predictor can be regarded as a solution to the variable selection problem when $\mathbb{P}(Y|X)$ is estimated accurately.

### 3.3  The Intersection Property and SDR

We have established that the Markov boundary is an optimal variable set for a learning algorithm under certain conditions. We now link Markov boundary with SDR theory. Note that Theorem 2.3.4 provides an important sufficient condition regarding the uniqueness of a Markov boundary. At the same time, Theorem 2.5.3 provides a different sufficient condition regarding the existence of a central DRS. We connect the two ideas by first proposing the following new definition, where $Z = (Y; X)$:

**Definition 3.3.1. (Linear intersection property).** The distribution $\mathbb{P}$ over variables $Z$ satisfies the linear intersection property if for every $\alpha_1, \alpha_2, \alpha_3, \alpha_4 \in \mathbb{R}^{(p+1) \times d}$ where $d \leq (p+1)$ such that $\alpha_1^T Z \perp\!\!\!\perp \alpha_2^T Z | (\alpha_3^T Z, \alpha_4^T Z)$ and $\alpha_1^T Z \perp\!\!\!\perp \alpha_4^T Z | (\alpha_3^T Z, \alpha_2^T Z)$, we also have $\alpha_1^T Z \perp\!\!\!\perp (\alpha_2^T Z, \alpha_4^T Z) | \alpha_3^T Z$. In other words, $\alpha_1^T Z \perp\!\!\!\perp \alpha_2^T Z | (\alpha_3^T Z, \alpha_4^T Z)$ and $\alpha_1^T Z \perp\!\!\!\perp \alpha_4^T Z | (\alpha_3^T Z, \alpha_2^T Z) \Rightarrow \alpha_1^T Z \perp\!\!\!\perp (\alpha_2^T Z, \alpha_4^T Z) | \alpha_3^T Z$. ∎

Note that the intersection property is an instance of the linear intersection property, when each $\alpha_i$ is a matrix whose diagonal elements are set to zero or one, and the off-diagonal elements are all zeros. Observe that for any $\alpha \in \mathbb{R}^{(p+1) \times d}$, if $\mathbb{P}(Z) > 0$, then $\mathbb{P}(\alpha^T Z) > 0$. Thus, if the distribution of $Z$ is strictly positive, then it satisfies the linear intersection property just like the original intersection property.

The following theorem guarantees the existence of a central DRS, if the linear intersection property holds:

**Theorem 3.3.2.** Let $S(\alpha)$ and $S(\phi)$ be DRSs for $Y|X$. If the joint distribution of $(Y; X)$ satisfies the linear intersection property, then $S(\alpha) \cap S(\phi)$ is also a DRS.

**Proof.** Note that a DRS must always exist, since $S(I_p)$ is a DRS. Let $S(\alpha)$ and $S(\phi)$ both be DRSs for $Y|X$. Let $\delta$ be a basis for $S(\alpha)$ and $S(\phi)$ so that $\alpha = (\alpha_1, \delta)$ and $\phi = (\phi_1, \delta)$. Now if $\alpha_1 = 0$ and $\phi_1 = 0$, then the conclusion is trivially true. Now let $\alpha_1 \neq 0$ and $\phi_1 \neq 0$. For notational convenience, define the following:

$$W = \begin{pmatrix} W_1 \\ W_2 \\ W_3 \end{pmatrix} = \begin{pmatrix} \alpha_1^T X \\ \phi_1^T X \\ \delta^T X \end{pmatrix}.$$

Since $S(\alpha)$ and $S(\phi)$ are both DRSs, we write the following relationship between the CDFs:

$$F_{Y|W} = F_{Y|W_1, W_3} = F_{Y|W_2, W_3}.$$

We need to also show that $F_{Y|W} = F_{Y|W_3}$. Note that the above relationships between the CDFs imply that $Y \perp\!\!\!\perp W | (W_1 \cup W_3)$ and $Y \perp\!\!\!\perp W | (W_2 \cup W_3)$ which in turn imply that $Y \perp\!\!\!\perp W_2 | (W_1 \cup W_3)$ and $Y \perp\!\!\!\perp W_1 | (W_2 \cup W_3)$. If the linear intersection property holds, we have $Y \perp\!\!\!\perp W_2 | (W_1 \cup W_3)$ and $Y \perp\!\!\!\perp W_1 | (W_2 \cup W_3)$ both imply that $Y \perp\!\!\!\perp (W_2 \cup W_1) | W_3$. This in turn implies that $Y \perp\!\!\!\perp W | W_3$ by Proposition 4.6 of (Cook 1998). Thus, $F_{Y|W} = F_{Y|W_3}$. ∎

### 3.4  Markov Boundary Discovery and SDR

Recall that a strictly positive distribution is reasonable to assume whenever there is uncertainty about the data. Thus, we assume that the joint distributions under consideration satisfy the linear intersection property from here on. This implies that $M$ is unique by Theorem 2.3.4. Without loss of generality (w.l.o.g.), we consider arranging $X$ such that $X = (M; X\backslash M)$. Let $p_M$ denote the number of variables in $M$. We can similarly consider arranging the rows of a matrix $\gamma \in \mathbb{R}^{p \times d}$ such that $\gamma = (\gamma_M; \gamma_{X\backslash M})$, where $\gamma_M$ has $p_M$ rows, and $\gamma_{X\backslash M}$ has $p - p_M$ rows. We claim the following:

**Theorem 3.4.1.** Suppose the joint probability distribution of $(Y; X)$ satisfies the linear intersection property. Then, there exists a central DRS $S(\gamma_M)$ for $Y|M$, and $S(\gamma)$ is the central DRS for $Y|X$, where $\gamma_{X\backslash M}$ is a matrix of all zeros.

**Proof.** Since the joint probability distribution of $(Y; X)$ satisfies the linear intersection property, then the distribution of $(Y; M)$ also satisfies the linear intersection property (as the property holds for any four linear combinations of variables), so we know that there exists a central DRS $S(\gamma_M)$ for $Y|M$ such that $Y \perp\!\!\!\perp M | \gamma_M^T M$ from Theorem 3.3.2. Together with the fact that $Y \perp\!\!\!\perp X | M$, we conclude $Y \perp\!\!\!\perp X | \gamma_M^T M$ from Proposition 4.6 of (Cook 1998) which implies that $Y \perp\!\!\!\perp X | \gamma^T X$, where $\gamma_{X\backslash M}$ is a matrix of zeros. Thus, $S(\gamma)$ is a DRS for $Y|X$. Now suppose there exists a matrix $\beta \in \mathbb{R}^{p \times d}$ such that $S(\beta)$ is also a DRS for $Y|X$ and $S(\beta) \subset S(\gamma)$. Then, $Y \perp\!\!\!\perp X | \beta^T X \Longrightarrow Y \perp\!\!\!\perp M | \beta^T X$. Note that w.l.o.g. we can arrange $\beta$ so that $\beta = (\beta_M; \beta_{X\backslash M})$. If $S(\beta) \subset S(\gamma)$, then $\beta_{X\backslash M}$ must be a matrix of all zeros. We then have the contradiction $Y \perp\!\!\!\perp M | \beta_M^T M$. Thus, we must have $S(\gamma) \subseteq S(\beta)$ for any DRS $S(\beta)$, so $S(\gamma)$ must be the central DRS for the conditional distribution of $Y|X$, which is guaranteed to exist by Theorem 3.3.2. ∎

Now that we know $S(\gamma)$ is the central DRS of $Y|X$, we next address the question of whether we can identify a unique $M$ from $\gamma$. We claim that the answer is yes:

**Theorem 3.4.2.** Suppose the joint probability distribution of $(Y; X)$ satisfies the linear intersection property, and let $S(\gamma)$ be the central DRS for $Y|X$ as defined in Theorem 3.4.1. Let $d'$ denote the number of column dimensions in $\gamma$. Then, we have $\sum_{i=1}^{d'} |\gamma_{j,i}| > 0$, when the row $j$ corresponds to a variable within $M$, and $\sum_{i=1}^{d'} |\gamma_{j,i}| = 0$, when the row $j$ corresponds to a variable within $X\backslash M$.

**Proof.** The linear intersection property holds, so a central DRS $S(\gamma)$ must always exist by Theorem 3.3.2. Also note that $S(\gamma) = S(\gamma_M)$ by Theorem 3.4.1. This implies that $\sum_{i=1}^{d'} |\gamma_{j,i}| = 0$, when the row $j$ corresponds to a variable within $X\backslash M$. Next, suppose there exists a row $j$ such that $\sum_{i=1}^{d'} |\gamma_{j,i}| = 0$, when the row corresponds to a variable within $M$. We then have the contradiction $Y \perp\!\!\!\perp X | \gamma_M^T M \Longrightarrow Y \perp\!\!\!\perp X | (M\backslash M_j)$, where $M_j \in M$. ∎

We thus find that variables in the Markov boundary are identified by discovering the central DRS and identifying any deviations from zero in the coefficients of $\gamma$.

### 3.5 *Markov Boundary Discovery with Ridge-Regularized Linear Models*

We now consider minimizing the following losses:

$$L_1(\alpha, \beta) = E\{u(\alpha + \beta^T X, Y)\}, \tag{E3}$$

$$L_2(\alpha, \beta) = E\{u(\alpha + \beta^T X, Y)\} + \lambda tr(\beta^T \Sigma_X \beta), \tag{E4}$$

where $u(\cdot,\cdot)$ is a convex functional, and $tr(\cdot)$ denotes the trace. Notice that the second loss has a covariance matrix $\Sigma_X$ added into the ridge regularization. The covariance matrix converts the ridge regularization into a new desired form whose purpose will become clear in the proof of Theorem 3.5.2.

Let $K$ denote the number the column dimensions of $\alpha$ and $\beta$, and let $k \in [1,2,\ldots,K] = \mathbb{K}$. We first provide the following analysis of (E3):

**Theorem 3.5.1.** Assume $\Sigma_X$ is positive definite, and let $S(\eta)$ be any DRS such that $E(X|\eta^T X)$ is a linear function of $\eta^T X$. If $(\alpha^*, \beta^*)$ minimizes (E3) and $\beta^*$ is unique, then $S(\beta_k^*) \subseteq S(\eta)$ for all $k \in \mathbb{K}$.

**Proof.** First, when $\eta = I_p$, then $E(X|\eta^T X)$ is always of linear function of $X$. Also, we have:

$$E[u(\alpha + \beta^T X, Y)] = E_{Y,\eta^T X} E_{X|Y,\eta^T X}[u(\alpha + \beta^T X, Y)]$$

$$= E_{Y,\eta^T X} E_{X|\eta^T X}[u(\alpha + \beta^T X, Y)],$$

where the last equality follows since $Y \perp\!\!\!\perp X | \eta^T X$. Recall that $u(\cdot,\cdot)$ is convex, so we apply Jensen's inequality as follows:

$$E_{Y,\eta^T X} E_{X|\eta^T X}[u(\alpha + \beta^T X, Y)] \geq E_{Y,\eta^T X}[u(\alpha + \beta^T E(X|\eta^T X), Y)]$$

W.l.o.g., we assume that $E(X) = 0$. Additionally, since $E(X|\eta^T X)$ is a linear function of $\eta^T X$, we apply Proposition 4.2 of (Cook 1998) so that:

$$L_1(\alpha, \beta) \geq L_1(\alpha, P_\eta(\Sigma_X)\beta),$$

where $P_\eta(\Sigma_X)$ is the projection matrix $\eta(\eta^T \Sigma_X \eta)^{-1} \eta^T \Sigma_X$. The conclusion now follows because the solution $\beta^*$ is unique. ∎

A similar argument is given in (Duan and Li 1991) and (Cook 1998), but when $E(X|\eta^T X)$ is a linear function of all $\eta^T X$ such that $Y \perp\!\!\!\perp X | \eta^T X$, and when $K = 1$.

Note that the assumption that $E(X|\eta^T X)$ is a linear function of $\eta^T X$ may not always hold. However, it usually holds when $p$ is large (Hall and Li 1993) and always holds when the distribution of $X$ is elliptically symmetric (Eaton 1986). The above theorem thus implies an important point: we can find a subset of the minimum DRS $S(\eta)$ such that $E(X|\eta^T X)$ is still a linear function using (E3), provided that $\beta^*$ is unique. Moreover, $\beta^*$ is guaranteed to include a subset of the Markov boundary if $E(X|\eta^T X)$ is a linear function when $S(\eta) = S_{Y|X}$ in the sense of Theorem 3.4.2.

We now claim that we can drop the assumption that $\beta^*$ is unique by adding ridge regularization with a covariance matrix as in (E4):

**Theorem 3.5.2.** Assume $\Sigma_X$ is positive definite, and let $S(\eta)$ be any DRS such that $E(X|\eta^T X)$ is a linear function of $\eta^T X$. If $(\alpha^*, \beta^*)$ minimizes (E4), then $S(\beta_k^*) \subseteq S(\eta)$ for all $k \in \mathbb{K}$.

**Proof.** The first part of the proof is similar to Theorem 3.5.1. We have:

$$E[u(\alpha + \beta^T X, Y)] = E_{Y,\eta^T X} E_{X|\eta^T X}[u(\alpha + \beta^T X, Y)].$$

Applying Jensen's inequality:

$$E_{Y,\eta^T X} E_{X|\eta^T X}[u(\alpha + \beta^T X, Y)] \geq E_{Y,\eta^T X}[u(\alpha + \beta^T E(X|\eta^T X), Y)]. \tag{E5}$$

Also consider, for any $k \in [1, 2, \ldots, K]$:

$$var(\beta_k^T X) = var\left(E(\beta_k^T X | \eta^T X)\right) + E\left(var(\beta_k^T X | \eta^T X)\right)$$
$$\geq var\left(E(\beta_k^T X | \eta^T X)\right). \tag{E6}$$

Putting (E5) and (E6) together, we have:

$$L_2(\alpha, \beta) \geq E_{Y, \eta^T X}[u(\alpha + \beta^T E(X | \eta^T X), Y)] + \lambda \sum_{k=1}^{K} var\left(\beta_k^T E(X | \eta^T X)\right) \tag{E7}$$
$$\Leftrightarrow L_2(\alpha, \beta) \geq E_{Y, \eta^T X}[u(\alpha + \beta^T E(X | \eta^T X), Y)] + \lambda tr\left(var(\beta^T E(X | \eta^T X))\right).$$

W.l.o.g., we assume that $E(X) = 0$. Additionally, since $E(X | \eta^T X)$ is a linear function of $\eta^T X$, we apply Proposition 4.2 of (Cook 1998) so that:

$$L_2(\alpha, \beta) \geq L_2(\alpha, P_\eta(\Sigma_X)\beta),$$

where $P_\eta(\Sigma_X) = \eta(\eta^T \Sigma_X \eta)^{-1} \eta^T \Sigma_X$ is the projection matrix.

Now let $S(\eta)$ be the minimum DRS such that $E(X | \eta^T X)$ is still a linear function of $\eta^T X$. If $S(\beta_k) \supset S(\eta)$, then $var\left(E(\beta_k^T X | \eta^T X)\right) > 0$. This implies that (E6) and (E7) are strict inequalities. As a result, such a $\beta_k$ cannot be used to minimize (E4). ∎

The above proof is a modification of Theorem 1 in (Li, Artemiou, and Li 2011) which was used in the context of binary SVMs. We believe the above theorem is powerful, since it applies to all convex losses regardless of the nonlinearities in the functional relationship.

The conclusions summarized in Section 3.1 follow from a straightforward synthesis of the above theorems. The linear intersection property ensures that a unique Markov boundary and a central DRS exist by Theorems 2.3.4 and 3.3.2, respectively. Since $\Sigma_X$ is positive definite, and $E(X | \eta^T X)$ is a linear function of $\eta^T X$ when $Y \perp\!\!\!\perp X | \eta^T X$, we can solve (E4) and find a $\beta^*$ such that $S(\beta_k^*) \subseteq S_{Y|X}$ for all $k \in \mathbb{K}$ by Theorem 3.5.2. We finally recover a subset of the parents, children and spouses of $Y$ from the non-zero coefficients of $\beta^*$ by Theorems 3.4.2 and 2.3.6, assuming $\beta^*$ is a non-zero matrix which is almost always the case in practice.

## 4 Experiments

### 4.1 Implementation

We now describe an implementation of the above ideas using the empirical version of (E4). Let $Y_n$ and $X_n$ denote two fixed $p$ by $n$ matrices of $n$ samples from the joint distribution of $(Y; X)$. The empirical version of (E4) can be written as follows:

$$\frac{1}{n} \sum_{i=1}^{n} \{u(\alpha + \beta^T X_{n,i}, Y_{n,i})\} + \lambda tr(\beta^T \hat{\Sigma}_X \beta). \tag{E8}$$

If the empirical covariance matrix $\hat{\Sigma}_X$ is non-singular, then (E8) can be minimized with standard packages for RRLMs such as Glmnet (Qian et al. 2013) or LIBSVM (Chang and Lin 2011) by standardizing $X_n$ so that $Z_n = \hat{\Sigma}_X^{-1/2}(X_n - \bar{X}_n)$. Then (E8) becomes:

$$\frac{1}{n}\sum_{i=1}^{n}\{u(\alpha + \gamma^T Z_{n,i}, Y_{n,i})\} + \lambda tr(\gamma^T \gamma), \tag{E9}$$

where $\gamma = \hat{\Sigma}_X^{1/2}\beta$. We then recover $\hat{\beta}^*$ by $\hat{\beta}^* = \hat{\Sigma}_X^{-1/2}\hat{\gamma}^*$.

We propose to solve (E9) and use permutation analysis to recover the p-values of each of the coefficients in $\hat{\beta}^*$. When the number of distinct samples is fewer than the number of dimensions $p$, we use the method of shrunken covariance estimators to estimate the covariance matrix (Ledoit and Wolf 2004, 2012). From here on, we use the term Covariance Ridge with Permutation (CRP) to refer to a permutation analysis by repetitively solving (E9) in order to recover coefficient specific p-values.

### 4.2 Toy Examples

We provide two toy simulations demonstrating the new theory. First, we obtained 1000 samples 100 times from 5 independent random variables with normal distributions with standard deviations that were drawn from a standard normal. We then created the random variable $Y$ using the following fifth order polynomial:

$$Y = \varepsilon + \theta_0 + \sum_{i=1}^{5}\theta_i X_i^i,$$

where the coefficients $\theta$ were also drawn from a standard normal, and $\varepsilon$ is a normally distributed error with small standard deviation 1E-4 such that $\varepsilon \perp\!\!\!\perp (X_1, \dots, X_5)$. We then added 1, 5, 15, 45, and 95 additional independent variables with normal distributions also with standard deviations drawn from a standard normal. The five variables $(X_1, \dots, X_5)$ thus comprise the Markov boundary of $Y$, and the additional variables are not part of the Markov boundary of $Y$.

CRP identified a subset of the five Markov boundary variables in all one hundred replicates across all numbers of additional variables. The mean number of correctly identified Markov boundary variables was 2.11 (SD: 0.751), and there was no significant decrease in the number of correctly identified Markov boundary variables as the number of additional variables increased. This is expected because the multivariate normal is elliptically symmetric, so $E(X|\eta^T X)$ is a linear function of any $\eta^T X$. Hence, the optimal $\beta^*$ obtained from minimizing (E4) is guaranteed to contain a subset of the Markov boundary from Theorem 3.5.2. For additional empirical evidence, we also ran a similar experiment, where the variables were drawn from a beta(2,2) distribution; recall that this distribution is also elliptically symmetric. We obtained similar results, except with a higher mean number of correctly identified variables at 4.47 (0.028).

In the next experiment, the explanatory variables were drawn from beta distributions, where the alpha and beta parameters were chosen by sampling from a uniform distribution between 0 and 5. Most instantiations of the beta distribution are not elliptically symmetric, but $E(X|\eta^T X)$ is almost always a linear function in the high dimensional setting, so we expect (E4) to still perform reasonably well even when we increase the number of additional variables. Nonetheless, we also expect that the method will not perform as well as with the beta(2,2) distributed variables. Indeed, we found that the mean number of correctly identified variables dropped to 4.18 (0.061) but again with no significant degradation as the number of additional variables increased.

### *4.3 Expert-Designed Models*

We evaluated CRP on datasets generated from four expert-designed discrete Bayesian networks (brief descriptions are given in Table 1). Note that many of the relationships between the variables in these networks are nonlinear. We compared the performances of the following three methods:

(1) The HITON-PC algorithm is a parent and child discovery method (Aliferis, Tsamardinos, and Statnikov 2003). Even though this method is not a Markov boundary discovery algorithm per se, it outperforms many existing algorithms on several metrics in identifying the true Markov boundary (Aliferis et al. 2010a, b). We used the Causal Explorer implementation (Statnikov et al. 2010) with the $G^2$ test. Additionally, we selected the $\alpha$ hyperparameter from the set {0.001, 0.01, 0.05, 0.10} using 5-fold cross-validation with an RBF-kernel SVM. The SVM's $C$ and $\sigma$ hyperparameters were in turn chosen within the folds from {0.1, 1, 10 100} and the median distance between samples multiplied by {2/3, 0.8, 1, 1.25, 1.5}, respectively. Finally, we set the $k$ hyperparameter for HITON-PC to 2, the largest value we could use so that the algorithm completed in a reasonable length of time in our experiments. In practice, values of $k$ up to 4 are recommended for HITON-PC (Aliferis et al. 2010a).

(2) The Kernel Backward Elimination (KBE) algorithm is, to our knowledge, the best performing Markov boundary *ranking* method (Strobl and Visweswaran 2013). Note that ranking may be more useful than selection when the number of selected variables is large. We used the default hyperparameter settings of $\sigma$ set to the median distance between samples and the ridge regularization set to 1E-4.

(3) The CRP algorithm was implemented by solving (E9) using the Glmnet package (Qian et al. 2013). We used a multinomial logistic loss or the MSE loss, when the former failed due to insufficient samples per group. Hyperparameters were determined by the cvglmnet function as provided in the package, and 1000 permutations were used to obtain approximate p-values.

In this experiment, our goal was to assess the ability of the above three algorithms in identifying a subset of the Markov boundary. Note that HITON-PC and KBE are guaranteed to discover the Markov boundary of any response in the infinite sample limit under their respective assumptions. However, CRP can discover only a subset of the Markov boundary of some unknown cardinality as noted in Theorem 3.5.2. We thus chose to evaluate the algorithms in their ability to discover a single variable in the Markov boundary. To do this, we randomly chose 50 variables as the responses from each Bayesian network, ran the three algorithms using sample sizes of {50, 100, 200, 300, 400, 500}, and then counted the number of times out of 50 that each algorithm correctly identified any Markov boundary variable. We repeated the experiments for each of the 4 networks and thus ran a total of 50×4=200 experiments per sample size. Recall that KBE and CRP output an ordering of the variables (the latter in terms of the smallest to largest p-values), so we only evaluated the $h$ lowest ranked variables for these two methods, where $h$ is the cardinality of the output from HITON-PC. We did not simply impose a p≤0.05 cut-off for CRP for reasons that will become clear in Section 4.5. The results of the experiments across multiple sample sizes are summarized in Figure 1.

| Network | Num. of Variables | Num. of Edges | Min / Max \|PC\| |
| --- | --- | --- | --- |
| Child10 | 200 | 257 | 1 / 8 |
| Alarm10 | 270 | 570 | 1 / 9 |
| Pigs | 441 | 592 | 1 / 41 |
| Gene | 801 | 972 | 0 / 11 |

**Table 1.** Descriptions of expert-designed discrete Bayesian networks. The fourth column gives the minimum and maximum cardinality of the parent and child set for each network.

To evaluate differences in performance regardless of sample size, we counted the number of times a Markov boundary variable was identified across all sample sizes and performed two-tailed Fisher's exact test. We used sample sizes of {50, 100, 200, 300, 400, 500} because these sample sizes are typical of the HDLSS datasets in biology to which RRLMs are applied. We found that HITON-PC and KBE outperformed CRP (both p<0.001) using the Bonferroni-adjusted threshold of 0.05/2=0.025. This result held even when we controlled for sample size and network using a $G^2$ test ($G^2$=81.67, p<0.001 with HITON-PC; $G^2$=39.51, p=0.004 with KBE). The results therefore suggest that CRP underperforms the nonparametric methods. However, examining the histograms in Figure 1 reveals that CRP is usually only outperformed by a small margin in the majority of sample sizes despite "only being linear."

### 4.4 Linear Gaussian Non-Recursive Structural Equation Models

Some real-world causal networks are likely to be cyclic or contain continuous variables. Examples include networks based on gene expression, where local causal discovery methods are most often applied. As a result, we evaluated CRP on datasets generated from linear Gaussian non-recursive structural equation models with independent errors (SEM-IEs), since our theoretical conclusions as summarized in Section 3.1 hold regardless of cyclicity. Moreover, it is known that linear Gaussian SEM-IEs and their stationary distributions (if they exist) satisfy the global Markov condition (Spirtes 1995). Recall also that the multivariate Gaussian distribution is strictly positive and therefore satisfies the intersection property by Proposition 2.2.2. The Markov boundary of the response in a stationary linear Gaussian SEM-IE thus uniquely includes the parents, children and spouses because these variables (but no proper subset of them) d-separate all other variables from the response.

We generated four networks using the default settings and default randomization methods in the TETRAD software package (Scheines et al. 1998). Moreover, unlike in Section 4.3, we equipped HITON-PC with the z-test (using Fisher's r-to-z transformation with the partial autocorrelation coefficient), and KBE with linear reproducing kernels due to the linearities in the network. Hyperparameters for HITON-PC were chosen by 5-fold cross validation using ridge regression with the best $\lambda$ obtained from the set {1E-6, 1E-4, 1E-2} as assessed within the folds. We also used the MSE loss for CRP. Both the task and the additional settings were otherwise kept the same as in Section 4.3. We therefore performed a total of 200 experiments per sample size. Note that the multivariate Gaussian distribution is elliptically symmetric, so we expect CRP to perform well, since its solution space is guaranteed to only contain a subset of the Markov boundary in the large sample limit by Theorem 3.5.2. Results are summarized in Figure 2 and indeed show that CRP's performance is competitive. CRP even outperforms HITON-PC (but not KBE) by a small margin by two-tailed Fisher's exact test (p=0.008 with HITON-PC; p=0.205 with KBE) and by $G^2$ test controlling for sample size and network ($G^2$=40.49, p=0.014 with HITON-PC; $G^2$=12.07, p=0.883 with KBE).

### 4.5 Application to Gene Expression Data

In the final set of experiments, we applied the algorithms to two real gene-expression datasets titled NOTCH1 and RELA (Statnikov et al. 2012). These datasets are valuable in that they come with gold standard solution sets consisting of experimentally verified genes that are downstream of the NOTCH1 and RELA transcription factors. The downstream genes were specifically identified using the following procedure. First, the samples were randomized to either a control or experimental group, where the transcription factor of interest was knocked down (e.g., by siRNA). A t-test was then performed with alpha level set to 0.05 to identify differentially expressed genes between the two groups. Second, the set of genes with promoter region-transcription factor binding was identified using genome-wide binding data (ChiP-on-chip for NOTCH1, and ChIP-sequencing for RELA). The final silver standard set of genes was obtained by overlapping the set of genes from the experimental knockdown data and the set of genes from the binding data.

Our goal was to compare the algorithms in their ability to identify the downstream genes of the transcription factors. We were also interested in the ease to which the algorithms' outputs could be interpreted by humans in order to help prioritize experimentation. We ran HITON-PC on the two datasets using the $G^2$ test (as suggested in (Aliferis et al. 2010a)), where each variable was first discretized into three bins. The $\alpha$ hyperparameter was chosen using 5-fold cross validation from {0.001, 0.01, 0.05, 0.10} using RBF kernel ridge regression with $\lambda$ from {1E-6, 1E-4, 1E-2} and $\sigma$ from the median distance between samples multiplied by {2/3, 0.8, 1, 1.25, 1.5}. HITON-PC ultimately provided an output of 1530 variables for NOTCH1 and 15375 variables for RELA. Similarly, CRP using an MSE loss and a p≤0.05 cut-off provided an output of 1045 variables for NOTCH1 and 921 variables for RELA. Interpreting such large sets of genes is clearly impractical for humans, and experimentally verifying them is even more impractical. Of course a simple heuristic solution would be to decrease the alpha threshold until the output contains a reasonable amount of variables. On the other hand, an alternative and perhaps a more principled strategy is to provide to the experimenter an ordered sequence of variables instead of an unordered set of variables such that the lowest ranked variables have the highest priority for experimental manipulation. We therefore compared KBE and CRP in detail, as both algorithms output variable rankings.

The KBE algorithm was run as described in Section 4.3, and the CRP algorithm was run with the MSE loss. Histograms of the rankings for the silver standard genes are shown in Figure 3. Note that the ideal output corresponds to a single probability mass at rank 1, and thus a more right skewed distribution denotes better performance. The histograms show that CRP outperforms KBE by a large margin, since the distributions of the rankings for CRP are more right skewed than those for KBE. Note that the histograms created from the output of KBE using a linear reproducing kernel were also left skewed.

## 5 Discussion

We have shown that the Markov boundary is an optimal solution to the feature selection problem for learners which must approximate the conditional distribution of $Y$ given $X$. However, empirical evidence suggests that identifying the Markov boundary in its entirety can be difficult for some problems involving high dimensional data. We can thus instead choose to identify a subset of the Markov boundary. By connecting Markov boundary and sufficient dimension reduction theory, we have shown that a modified form of RRLMs can get very close to identifying a subset of the Markov boundary. This fact holds under a variety of losses and nonlinearities in the functional relationships between the response and explanatory variables. In practice, the modified RRLMs compare favorably to state-of-the-art full Markov boundary discovery methods in both performance and interpretability on an experimentally verifiable task with two gene expression datasets.

The CRP algorithm may have outperformed KBE on gene expression data due to a different reason besides attempting to solve an easier problem. Specifically, ranking variables by their p-values may in practice be superior than ranking them by their residual errors (when removed). This preliminary claim is partially supported by the observation that CRP also outperforms KBE with linear reproducing kernels which is in fact equivalent to linear ridge regression (Fukumizu and Leng 2014, Fukumizu, Bach, and Jordan 2009). Naturally, one may next wonder whether combining p-values with backward elimination can yield even better performance. Unfortunately, obtaining p-values at every iteration with a backward elimination method is often too computationally expensive in practice, so this strategy is not tractable unless the asymptotic distributions of the coefficients can be predefined under the null.

Besides its advantages, CRP also has several limitations. First, we must emphasize that this method does not replace Markov boundary discovery methods which are guaranteed to discover the Markov boundary in its entirety in the infinite sample limit (under their respective assumptions). Nonetheless, we hope that the ideas presented in this paper suggest that RRLMs are more useful than previously demonstrated. Second, we were unable to run the method on high-dimensional datasets with several tens of thousands of variables on desktop computers due to the memory requirements of the covariance matrix. This limits the applicability of the method, and we are not currently aware of ways to alleviate the issue. Third, the application of the method to gene expression data may seem premature given current understanding about d-separation in cyclical causal models. Indeed, the global Markov condition is only known to apply to the linear Gaussian as well as to discrete SEM-IEs when cyclicity exists (Spirtes 1995, Pearl and Dechter 1996). However, the method appears to perform well in practice, either because the global Markov condition holds in the tested cases or due to other reasons. For example, Spirtes defined another separation condition which holds more generally with non-linear SEM-IEs via d-separation in *collapsed graphs*, which we now call *cd-separation* for collapsed d-separation (Spirtes 1995). Unfortunately, cd-separation may not be the finest separation condition across all stationary distributions, in the sense that there are separation conditions which may entail additional conditional independencies; indeed, cd-separation implies d-separation but the converse is not true, so cd-separation is coarser than d-separation. Equivalently, we can state that d-connection implies cd-connection. We believe that the modified RRLMs may, in the worst case, be discovering a subset of those variables which are cd-connected with the response. As a result, the modified RRLMs are also correctly discovering some of those variables which are d-connected with the response.

A variety of interesting open questions remain. We do not know if the same conclusions hold if the covariance matrix is dropped from the ridge, or if the lasso is used instead. Moreover, it may be interesting to know whether certain losses have stronger theoretical guarantees than others in identifying a larger portion of the Markov boundary. Nonlinear methods such as those based on reproducing kernels may further allow us to drop some of the assumptions in Theorem 3.5.2. Finally, it may be worthwhile to explore violations of the intersection property, where a central DRS may not exist. We conjecture that the solution to a modified RRLM will be located within the union of all minimal DRSs and thus can be used to identify a subset of the union of all Markov boundaries.

In conclusion, we have presented theoretical and experimental results concerning RRLMs that strengthen their interpretation as a causal discovery method. Although modified RRLMs are not guaranteed to detect a subset of the Markov boundary in the large sample limit, they almost do so in theory and often do so in practice.

# 6 Acknowledgments

Research reported in this publication was supported by the National Institutes of Health under award number U54HG008540 as well as the National Library of Medicine under award number 5T15LM007059-28. The content is solely the responsibility of the authors and does not necessarily represent the official views of the National Institutes of Health or the National Library of Medicine.

## 8 Figures

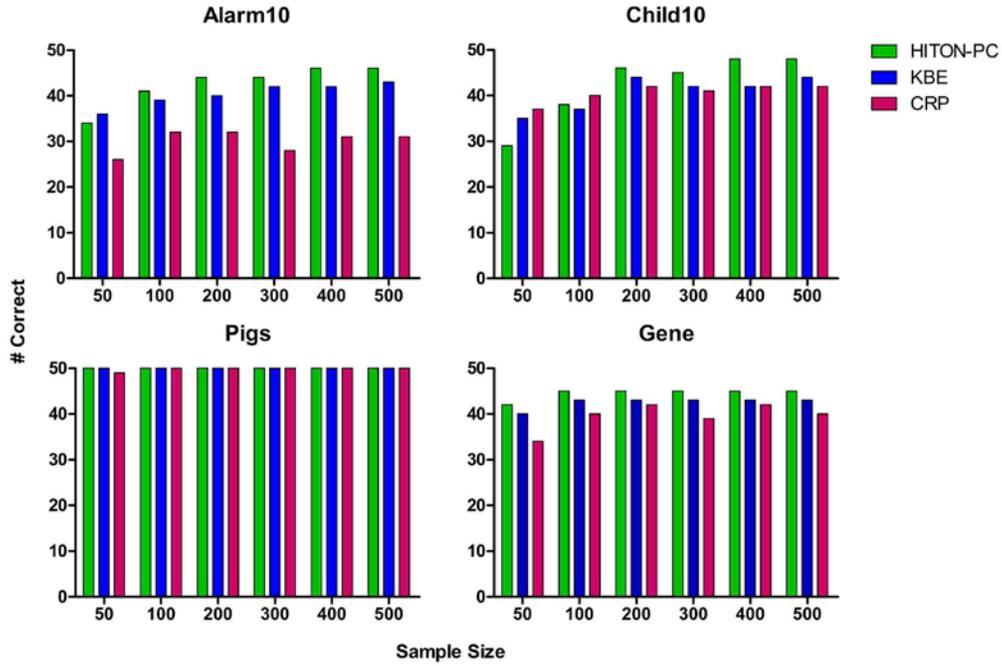

**Figure 1.** The number of times HITON-PC, KBE, and CRP correctly detected a single variable in 50 Markov boundaries across four expert-designed discrete Bayesian networks. HITON-PC and KBE outperform CRP by a small margin across most sample sizes.

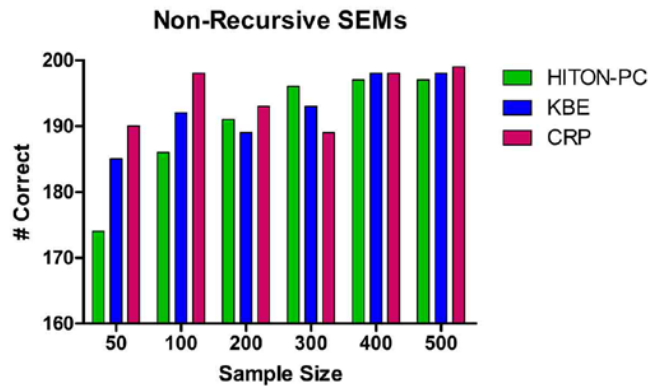

**Figure 2.** Same as Figure 1 except the 50 Markov boundaries were obtained from four linear Gaussian non-recursive SEM-IEs. The results are combined across the four SEM-IEs. CRP outperforms HITON-PC by a small margin in detecting a single variable within the Markov boundary.

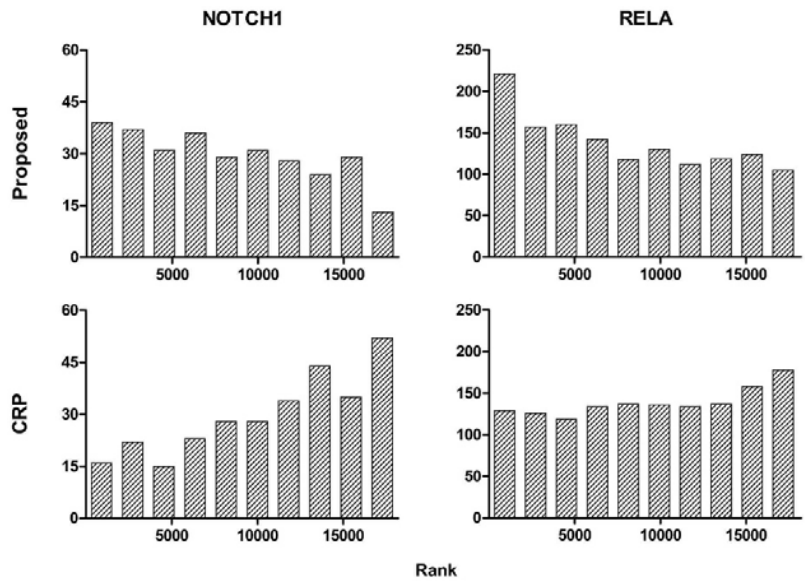

**Figure 3.** Histograms of the rankings obtained from CRP (top row) and KBE (bottom row) for the NOTCH1 (first column) and RELA (second column) datasets. CRP outperforms KBE, since the histograms for CRP are more right skewed than those for KBE.